\documentclass[conference]{IEEEtran}
\IEEEoverridecommandlockouts
\usepackage{cite}
\usepackage{amsmath,amssymb,amsfonts}
\usepackage{algorithmic}
\usepackage{float}
\usepackage{graphicx}
\usepackage{textcomp}
\usepackage{xcolor}
\newtheorem{theorem}{Theorem}[section]

\def\BibTeX{{\rm B\kern-.05em{\sci\kern-.025em b}\kern-.08em
    T\kern-.1667em\lower.7ex\hbox{E}\kern-.125emX}}

\begin{document}

\title{Geometric Control of a Robot's Tool\\
}

\author{\IEEEauthorblockN{Bousclet Anis and Belhadjoudja Mohamed Camil}
\IEEEauthorblockA{\textit{Control Engineering Department} \\
\textit{National Polytechnic School}\\
\textit{Students in Master of Engineering} \\
\textbf{Process Control Laboratory}\\
Algeria, Algiers \\ 
anis.bousclet@g.enp.edu.dz, mohamed\_camil.belhadjoudja@g.enp.edu.dz }
}
\maketitle

\begin{abstract}
The goal of this paper is to present a rigorous and intrinsic formulation of a riemannian PD-regulator of the robot's tool, The first one is based upon the Lasalle's invariance principle, we use it to control the tool's position in the workspace under the assumption of absence of singularities in configuration space, The second method deals with geometrical constraints on the trajectory of the robot's tool with the same assumption, we construct a unique orthogonal force that is viewed as a gravitational force that keeps the tool constrained, We also present a variation of the first method in the case of double pendulum based on the Lyapunov stability theorem. With this modification, we control the tool and the difference between the two angles, we did simulations on a two-link manipulator that shows the efficiency of the presented methods.
\end{abstract}


\section{Introduction}
Robotic manipulators are very useful in industries, agriculture, medicine and other important domains, they allow us to gain time, precision and efficacy, the main problem of control robotic systems is the non-linearity of the dynamics, there are some books and papers that model the configuration space of the robot as an Euclidean space and apply the principle of least action to give the Euler-Lagrange Equations [9] [10] [11] [12] [13] [14] [15] [16] [17], and apply Lyapunov second method to control the robot with Euclidean PD-Regulator, in the other hand Arnold, Abraham and Marsden showed the importance of Riemannian geometry in classical mechanics [5] [6], this was followed by some works of Selig and Murray's team for geometric modeling of robotic systems [20] [21], recently, the team of Suguru Arimoto presented very interesting results based upon Riemannian geometry [22] [23] [24] [25], the wonderful work of Bullo and Lewis was to give a rigorous and intrinsic formulation and proofs of the results of Arimoto's team, excepted some results concerning the control of the tool or end-point of the robot, and the control ensuring the geometric constraint of the tool, here we present an intrinsic and rigorous formulation of these results, giving proofs using tools of Riemannian geometry [4] [8] [9] and dynamical systems [3] [11], we also connect the absence of singularity [7] of the tool function with the validation of these methods.\\
\section{Notations and preliminaries}
\subsection{Geometric modeling of robotic systems}
A robot is a set of $s$ solids connected by joints. A priori, it may seem that we need $6s$ parameters to describe the evolution of this system, however this set is subject to holonomic constraints expressed in terms of $r$ submersion equations [5] [6], which reduces the number of coordinates needed from $6s$ to say $n$ variables [7]. as a simple pendulum is constrained to have a constant distance from the origin gives a circle as configuration space, we can prove under some assumptions [1] [7] that the configuration space is a manifold $M\subset (\mathbf{R}^3\times SO_{3}(\mathbf{R}))^s$ of dimension $n=6s-r$ for which the usual coordinates in analytical mechanics are local coordinates given from the chart of the manifold. This dimension is called the robot's degree of freedom DOF and is the exact number of coordinates we need to give a description of the robot, each value of the variable is called a configuration of the robot. The phase space is the space of the initial positions and speeds $(q,v)$ of a system, it is the tangent bundle $TM$ of the configuration space, we refer to [1] [7] for more details. \subsection{Robot's Tool function, Singularities }
The terminal organ of a robot is what is called the tool. By knowing the rotation matrix and the position of the center of gravity of each solid, it is possible to determine the tool function $x : M \to \mathbf{R}^3$ which, at each configuration of the system, gives the position of the terminal organ. The workspace is defined as the set of all the points that the tool can reach, we denote it by $E_{T}:=x(M)\subset \mathbb{R}^3$.\\
We call a singular point [1] [7] of a robot any point $q\in M$ such that $dx_{q} : T_{q}M \to \mathbf{R}^3$ is not onto, when we move in a neighborhood of $q$ the tool does not move in all possible directions. A singular value is a point $x_{d} \in E_{T}$ such that $\exists q\in x^{-1}(\lbrace x_{d}\rbrace )$ such that  $dx_{q} : T_{q}M \to \mathbf{R}^3$ is not onto, which means that there is at least one configuration in the preimage of $\lbrace x_{d}\rbrace$ that is a singular point. A regular point is any point $q\in M$ such that $dx_{q} : T_{q}M \to \mathbf{R}^3$ is onto. A regular value is a point $x_{d} \in E_{T}$ such that $\forall q\in x^{-1}(\lbrace x_{d}\rbrace )$, $dx_{q} : T_{q}M \to \mathbf{R}^3$ is onto, which means that every configuration in the preimage of $\lbrace x_{d}\rbrace $ is a regular point.
\subsection{Riemannian structure of the configuration space }
Let $(M,g)$ be a Riemannian manifold, and $\gamma : [t_{0},t_{1}]\subset \mathbf{R}\to M$ a smooth curve. We call kinetic energy [1] [5] [6] [19] of $\gamma$ the functional: 
$$E(\gamma )=\frac{1}{2}\int_{t_{0}}^{t_{1}}g_{\gamma (t) }(\gamma '(t), \gamma '(t))dt$$
In our work, we see the configuration space as being a Riemannian manifold whose metric is the kinetic energy that we know from rigid body dynamics in order to understand the trajectory of the systems in absence of external forces as geodesics of the kinetic energy. In order to compute this metric, we add the kinetic energies of rotation and translation of each solid. Let $m_{j}$ be the mass of the $j^{th}$ solid, $x_{Gj}$ its center of mass, $I_{Gj}$ its tensor of inertia about $x_{Gj}$, $R_{j}$ its rotation matrix in a fixed frame and $\Omega  (R_{j}(t),R'_{j}(t))$ its instantaneous rotation vector. We will write $\Omega_{j} '$ instead of $\Omega '(R_{j}(t),R'_{j}(t))$.  The formula for the metric is then [1] [6] [12] [15] [17] [19]: 
$$ g_{\gamma (t)}(\gamma '(t), \gamma '(t)) =  \sum_{j=1}^{s}\lbrace m_{j}[x_{Gj}'(t)]^2 + \Omega_{j} 'I_{Gj}\Omega _{j}'\rbrace $$
In addition to the kinetic energy, we will choose a smooth function $U:M\to R$ called potential energy. \\
\subsection{Equation of motion }
 Let $[t_{0},t_{1}]$ and $]-\epsilon , \epsilon [$ be subsets of $\mathbf{R}$. For any smooth curve $\gamma : [t_{0},t_{1}] \to M$, we define the action of this curve by the formula : 
$$ S(\gamma ) = \int_{t_{0}}^{t_{1}}[\frac{1}{2}g_{\gamma (t)}(\gamma '(t), \gamma '(t)) - U(\gamma (t))]dt $$
 According to the principle of least action [10], the robot will evolve in such a way that the curve $\gamma  : [t_{0},t_{1}] \to M$ which describes the evolution of the robot's variables on the configuration space would minimize the action $S$.\\ \\ If we assume that the system evolves from the point $p_{1}\in M$ to the point $p_{2}\in M$ in an interval of time $[t_{0},t_{1}]$, then $\gamma : [t_{0},t_{1}] \to M$ must satisfy [1] [5] [6] [8] [9] [10] [19]: 
$$S(\gamma ) = \inf_{\alpha \in \Omega (p_{1},p_{2})} S(\alpha )$$
 Where $\Omega (p_{1},p_{2})$ is the set of smooth maps $\alpha : [t_{0},t_{1}] \to M$ such that $\alpha (t_{0})=p_{1}$ and $\alpha (t_{1})=p_{2}$. \\ \\
Using some calculus of variations [1] [5] [6] [8] [9] [10] [19], we find the following equation that is a generalized equation of geodesics and of the newton equation that says how the potential and the curvature of the configuration space affects the trajectory of the robot:
$$\frac{D\gamma'(t)}{Dt}=-grad_{g}(U)(\gamma(t))$$
We have the following theorem which assures us the existence and the uniqueness of a solution to the generalized equation of geodesics given certain initial conditions [3] [4] [7].\\
\begin{theorem}{Cauchy-Lipschitz}\\
For every initial conditions $(q,v)\in TM$, there exists a unique maximal curve defined on an open interval $I_{v}\subset \mathbf{R}$ containing $0$. This curve $\gamma : I_{v} \to M$ starts in $x$ with initial speed $v$ and satisfies the equation (1). 
\end{theorem}
\section{Robot Control}
To control the robot, we act on its acceleration by adding a control term [1] [12-26], a fictious control law is a smooth map from $I$ to $TM$ such that for $t\in I$, $u(t)\in T_{\gamma(t)}M$
$$\frac{D\gamma'}{Dt}(t)=-\nabla U(\gamma(t))+u(t)$$
truly speaking [19], the real control is an application $u:I\rightarrow{T^{*}M}$ such that $u(t)\in T_{\gamma(t)}M^{*}$ for all $t\in I$, and it acts on the acceleration of the robot by $\tilde{u}(t)=g_{\gamma(t)}^{\sharp}(u(t))$ which is an application from $I$ to $TM$ such that for all $t\in I$ $\tilde{u}(t)\in T_{\gamma(t)}M$.  \\
\subsection{Kinetic energy theorem }
The quantity $E(t) = (1/2)g_{\gamma (t)}(\gamma '(t), \gamma '(t)) + U(\gamma (t)) $, with $\gamma $ a solution of the equation of motion is called total energy of the robot. We have the following results [1] [19]: \\ \\
\begin{theorem}
\textit{Let $\gamma : I\subset \mathbf{R}\to M$ be a solution of (5), then $\forall t_{1},t_{2}\in I$ :}
$$E(t_{1})-E(t_{0}) = \int_{t_{0}}^{t_{1}} g_{\gamma (t)}(u(t),\gamma '(t))dt$$ 
\textit{Or locally :}
$$\frac{d}{dt}[E(t)]=g_{\gamma (t)}(u(t), \gamma '(t))$$\end{theorem}
This theorem gives us insights about how the total energy of a robot varies with the control law, and this is very important for the study of stability. 
\subsection{Control of tool's position }
In this control method [1] [24] [25], the reference is the tool's position $x_{d}$ in the workspace. We want to make $\lim_{t\to \infty} x({q}(t)) = x_{d}$ with zero configuration velocity, where $x : M \to E_{T}$ is the tool function and $q : I\subset \mathbf{R} \to M$ is the configuration variable. Assume that there are no singular points. The idea here is to make each point of $x^{-1}(\lbrace x_{d}\rbrace )$ an equilibrium point. We achieve this goal by choosing a control law that compensates the conservative forces, dig holes in the neighborhood of configurations for which the tool is in the desired position and stabilize those equilibrium points. In order to dissipate the energy, we add friction forces. \\ \\ The control law is then : $$u(q,v) = \nabla U(q)-k.v-\nabla V(q)$$ Where $$V(q)=\frac{k_{1}}{2} || x(q)-x_{d}||^2, \ \ k_{1}>0$$ The closed-loop dynamics becomes : 
$$\frac{D\gamma'}{Dt}(t)=-grad_{g}(V)(\gamma(t))-k.\gamma'(t)$$
The equilibrium points are exactly the points of $x^{-1}(\lbrace x_{d}\rbrace )$, we have $grad_{g}(V)(q)=0$ if and only if $dV_{q}(v)=0$ for all $v\in T_{q}M$, using the fact that $dV_{q}(v)=k_{1}<x(q)-x_{d},dx_{q}(v)>$ we see that $grad_{g}(V)(q)=0$ is equivalent to $q\in x^{-1}(\left\{x_{d}\right\})$. In terms of energy, we have : $$\frac{d}{dt}[\frac{1}{2} g_{\gamma(t) }(\gamma '(t),\gamma '(t))+V(\gamma (t))] =-k.|\gamma'(t)|_{g}^{2}$$ Because of the fact that the energy is a proper function on the tangent bundle, that it is decreasing in all trajectories of the robot, and the largest invariant subset in the Lasalle's invariance principle $\Omega = x^{-1}(\lbrace x_{d} \rbrace ) \times \lbrace 0 \rbrace $, we find from here that the configuration converges to $x^{-1}(\lbrace x_{d} \rbrace )$ and the velocity to $0$. This means that $x(q(t))$ converges to $x_{d}$ by continuity of the tool function.\\
so we have the folowing result :\\\\
\fcolorbox{black}{lightgray}{
\begin{minipage}{0.48\textwidth}

\begin{theorem}
let $(M,g)$ be a compact riemannian manifold that models a robot, and the tool function $x:M\rightarrow{\mathbb{R}^{3}}$ is without singularities, so for all $x_{d}\in E_{T}$ the feed-back control law $u(q,v)=grad_{g}(U)(q)-k.v-grad_{g}(V)(q)$ for $k>0$ make the tool reach $x_{d}$ with zero configuration velocity. 
\end{theorem}
\end{minipage}}
\subsection{Control with geometrical constraint on the tool's position :} From here on we suppose that $U=0$ ( after compensating it with the control law we have $\frac{D\gamma'}{Dt}=u$) and we want the tool to stay in an orientable surface $S\subset E_{T}$, let $\Phi:\mathbb{R}^{3}\rightarrow{\mathbb{R}}$ a submersion in $E_{T}$ and $S=\Phi^{-1}(\left\{0\right\})$.  \\
Let $x_{d} \in S$, we suppose that there is no singular point, this ensure that $\Psi=\Phi\circ x$ is a submersion in $M$, and $N=\Psi^{-1}(\left\{0\right\})=x^{-1}(S)\subset M$ will be a Riemannian sub-manifold of $M$.\\ 
so we have $x:N\subset M\rightarrow{S\subset E_{T}\subset\mathbb{R}^{3}}$, this will give $dx_{q}:T_{q}N\rightarrow{T_{x(q)}S}$ for all $q\in N$, it says that if we want the velocity of the tool be tangent to $S$, the configuration must be tangent to $N$. \\
The tool's displacement constraint on the surface has become a constraint on the configuration which must stay in $N$.\\
\subsubsection{The contact force}
The main idea [1] [24] [25] is to prove uniqueness and existence of the normal component of a feedback control law such that if the tool is initialized in $S$ with a tangent velocity, it will remain in $S$ all the time, suppose we have such a control $u$, so $\gamma:I\rightarrow{N}$, we have :
$$g_{\gamma}(grad_{g}(\Psi),\frac{D\gamma'}{Dt})=g_{\gamma}(u,grad_{g}(\Psi))$$ using the Leibnitz rule and the fact that for all $q\in N$ $grad_{g}(\Psi)(q).\mathbb{R}=T_{q}N^{\perp}$ we have : $$g_{\gamma}(u,grad_{g}(\Psi(\gamma))=-g_{\gamma}(\frac{Dgrad_{g}(\Psi)(\gamma)}{Dt},\gamma')$$ and this give uniqueness of the orthogonal component of the feedback law, adding a tangent feedback law to ensure the same conditions as in the unconstrained problem, we have $u=\lambda.grad_{g}(\Psi)+u_{//}$ with $\lambda= \frac{-g_{\gamma}(\frac{Dgrad_{g}(\Psi)(\gamma)}{Dt},\gamma')}{|grad_{g}(\Psi)(\gamma)|_{g}^{2}}$, and $u_{//}(\gamma,\gamma')\in (grad_{g}(\Psi)(\gamma).\mathbb{R})^{\perp}$.
$$\frac{D^{N}\gamma '}{Dt}=u_{//}(\gamma , \gamma ')$$
With $u_{//}(\gamma, \gamma') \in T_{\gamma}N$. \\ \\
on the other hand we have : $$grad(U)(\gamma)_{//}=P_{T_{\gamma}N}(grad_{g}(U)(\gamma))=grad_{g^{N}}(U|_{N})(\gamma)$$ With this remark, we can apply all our previous results just by replacing in the hypothesis $M$ by $N$, because we have eliminated the constrained with the feedback orthogonal control law.\\ \\ 
\fcolorbox{black}{lightgray}{
\begin{minipage}{0.48\textwidth}

\begin{theorem}
let $(M,g)$ be a compact riemannian manifold that models a robot system, let $S\subset E_{T}$ be a surface such that there exists $\Phi:\mathbb{R}^{3}\rightarrow{\mathbb{R}}$ submersive in $E_{T}$ such that $S=\Phi^{-1}(\left\{0\right\}$, and that there is no singularities of the tool function $x:M\rightarrow{\mathbb{R}^{3}}$, so there exists a unique normal component of the feed-back control law such that if the tool is initialized in $S$ with tangent velocity, it will remain on $S$ for all time, explicitely we have :
$$u^{\perp}(q,v)=\frac{-g_{q}(D_{v}grad(\Psi)(q),v)}{|grad_{g}(\Psi)(q)|_{g}^{2}}grad_{g}(\Psi)(q)$$ and all control law of the form $u(q,v)=u^{\perp}(q,v)+u_{//}(q,v)$ such that $u_{//}(q,v)\in (grad(\Psi)(q).\mathbb{R})^{\perp}$ will accomplish the same task.
\end{theorem}
\end{minipage}}\\ \\
\textbf{proof :}\\
let $q_{0}\in N$ and $v\in T_{q_{0}}N$, we consider the real function of real variable $f:I\rightarrow{\mathbb{R}}$ by $f(t)=\Psi(q(t))$, the initial conditions gives $f(0)=f'(0)=0$, we compute $f''(t)$, using the fact that $f'(t)=g_{q(t)}(grad(\Psi)(q(t)),q'(t))$ and the Leibniz rule, also the equation of motion of the robot and the formula of $u$, we have $f''(t)=g_{q(t)}(\frac{Dgrad(\Psi)(q)}{Dt}(t),q'(t))+g_{q(t)}(grad(\Psi)(q(t)),\frac{Dq'}{Dt}(t))$ replacing $\frac{Dq'}{Dt}=u^{\perp}+u_{//}$ we find $f''(t)=0$, this finishes the proof. 
\\
\subsubsection{Control of tool's position under constraint}
We denote $v_{//}$ by $P_{(\nabla \Psi(q))^{\perp}}(v)$ for $v\in T_{q}M$, using the feedback control law $$u^{//}(q,v)=-k.v_{//}-grad_{g}(V)(q)_{//}$$ we can make the tool converge into $x_{d}\in S$ while remaining in $S$ for all times, the only point that is not clear is when a system is constrained, it will have more equilibrium positions, under a geometric hypothesis on the constrained surface $S$. \\ \\
\fcolorbox{black}{lightgray}{
\begin{minipage}{0.48\textwidth}

\begin{theorem}
suppose that the tool function is without singularities and suppose also that for all $q\in N$ such that $x(q)-x_{d}\in T_{x(q)}S^{\perp}$ we have $x(q)=x_{d}$, so the feed-back control law $u=u^{\perp}+u^{//}$ $$u^{//}(q,v)=-k.v_{//}-grad_{g}(V)(q)_{//}$$ will make the robot's tool initialized in $S$ with tangent velocity converging to $x_{d}\in S$ remaining in $S$ for all time. 
\end{theorem}
\end{minipage}}\\ \\
\textbf{Proof :}\\
We conclut by theorem 3.3 that the tool remains in $S$, and so $\gamma:I\rightarrow{N}$, by the orthogonal projection of the equation of motion on $\nabla\Psi(\gamma(t))=T_{\gamma(t)}N$ we have 
$$\frac{D^{N}\gamma'}{Dt}=-P_{T_{\gamma(t)N}}(\nabla V(\gamma(t)))-k.\gamma'(t).$$ The absence of singularities let us conclude that the assumption about $S$ is equivalent to the fact that "for $q\in N$ that satisfies $<dx_{q}(v),x(q)-x_{d}>=0$ for each $v\in T_{q}N$, so $q\in x^{-1}(x_{d})$", and this clearly means that critical points of $V$ in $N$ are $x^{-1}(x_{d})$, the Lasalle's invariance principle concludes. 
\section{Simulation examples}
\subsection{The two link manipulators}
\subsubsection{General presentation}
\begin{figure}[H]
\centerline{\includegraphics[width=90mm]{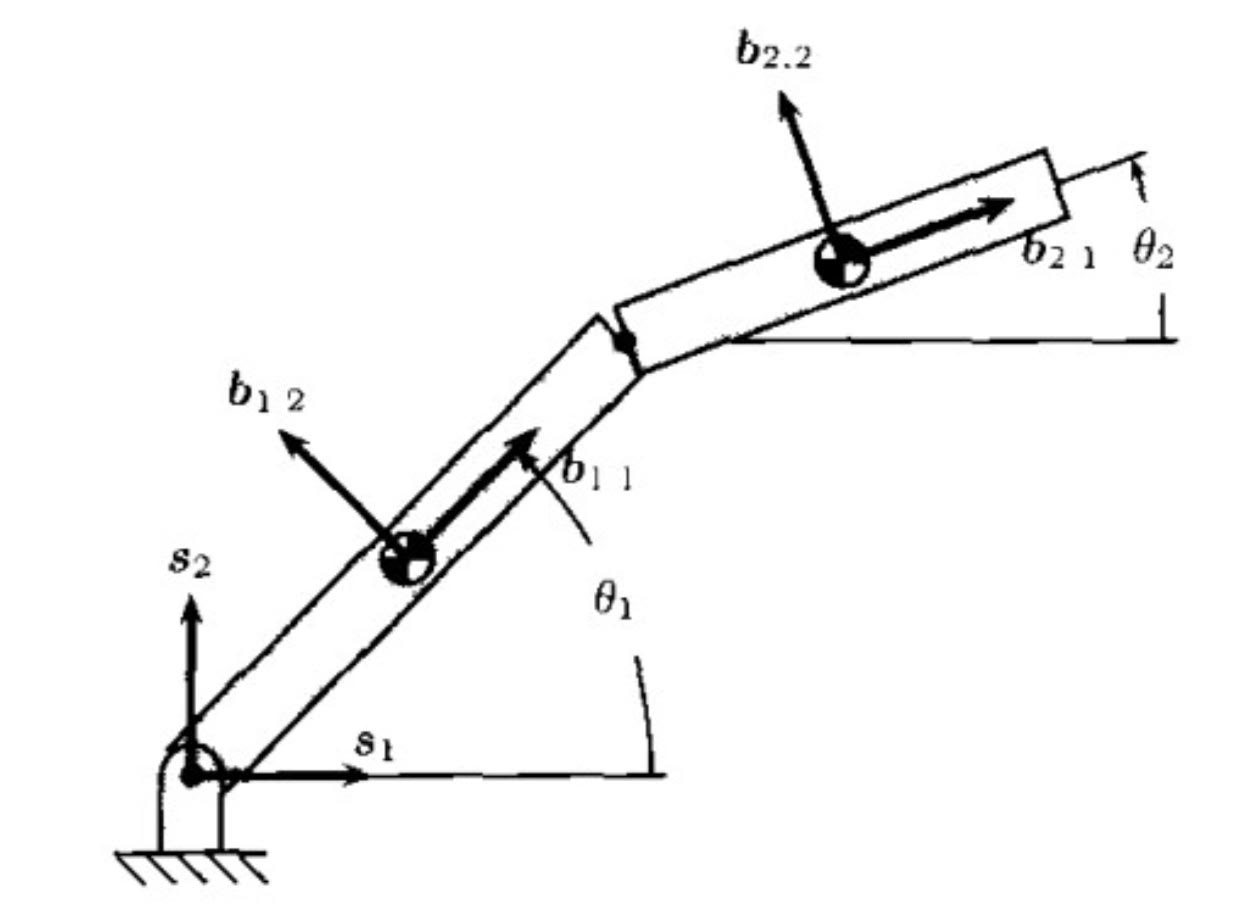}}
\caption{The two-link manipulator.}
\label{fig}
\end{figure} 
The two-link manipulator [19] is an articulated robotic arm made up of two solid links and two rotary joints. The first link has a mass $m_{1}$, a center of mass $x_{G1}$, a moment of inertia $J_{1}$ with respect to $x_{G1}$ and a length $l_{1}$. The second one has a mass $m_{2}$, a center of mass $x_{G2}$, a moment of inertia $J_{2}$ with respect to $x_{G2}$ and a length $l_{2}$. The end of the first link moves along the circle with radius $l_{1}$ and whose center is the first joint. The end of the second link moves along the circle with radius $l_{2}$ and whose center is the second joint (we will sometimes use link rather than end of the link, the meaning will be clear depending on the context of the sentence). The robot’s configuration is determined by the position of each link on the corresponding circle. \\

\subsubsection{Configuration space of the two-link manipulator}
The configuration space of the two-link manipulator is the Torus $T^{2}=\mathbb{S}^{1}\times \mathbb{S}^{1}$ [19], this manifold is a priori a sub-manifold of $\mathbb{R}^{4}$ but it can be embedded in $\mathbb{R}^{3}$ to give the usual "torus" that is known as a donut [7], In order to facilitate the computations and make them graphically more understandable, we will use local coordinates. The Torus being a manifold of dimension $2$, Instead of using a point of the Torus $q = [(x_{1}, y_{1}), (x_{2}, y_{2})]$, it will be better for us to use two real numbers which we will denote by $\theta_{1}$ and $\theta_{2}$. \\  The numbers $\theta_{1}$ and $\theta_{2}$ will be determined by the map $\varphi : T^{2} \to \mathbf{R^2}$ defined by [19]:  
$$\varphi (q) = (atan(x_{1},y_{1}),atan(x_{2},y_{2})) = (\theta_{1},\theta_{2})$$
Graphically, $\theta_{1}$ is the angle between the first link and the horizontal and $\theta_{2}$ the angle between the second link and the horizontal. These angles are counted positive in the counterclockwise direction.\\ 
\begin{figure}[H]
\centerline{\includegraphics[width=90mm]{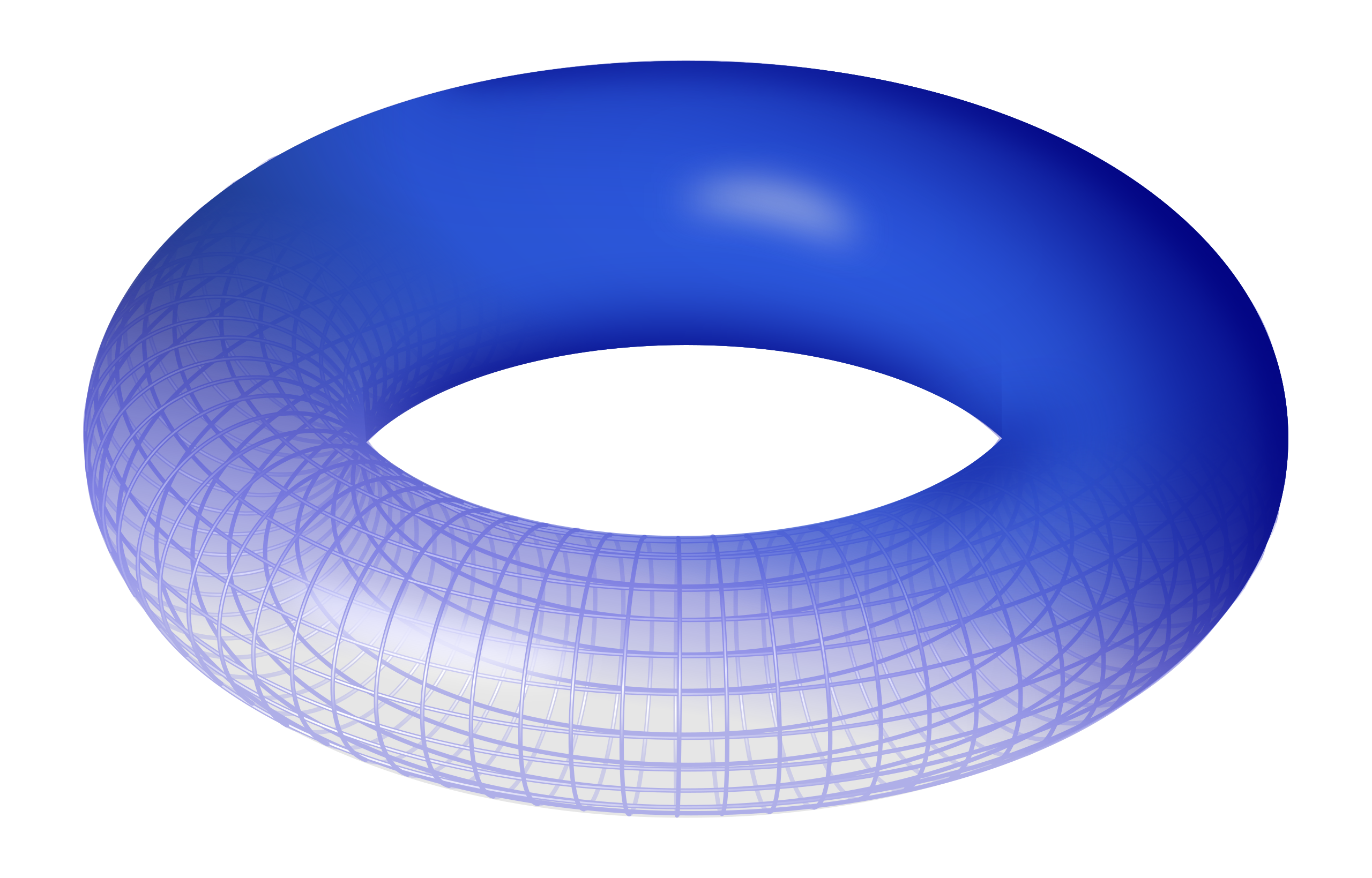}}
\caption{The configuration space $T^{2}$.}
\label{fig}
\end{figure}
\subsubsection{Riemannian structure of the torus}
We start by calculating the kinetic energy of the robot. We can show that it is given by the formula [19]: 
\begin{IEEEeqnarray}{rCl}
E & = & \frac{1}{8}(m_{1}+4m_{2})l_{1}^{2}\dot \theta_{1}^2 \nonumber \\ 
&&+\> \frac{1}{8}m_{2}l_{2}^2\dot \theta_{2}^2 \nonumber\\
&&+\> \frac{1}{2}m_{2}l_{1}l_{2}cos(\theta_{1}-\theta_{2})\dot \theta_{1}\dot \theta_{2} \nonumber\\ 
&&+\> \frac{1}{2}J_{1}\dot \theta_{1} +  \frac{1}{2}J_{2}\dot \theta_{2}\nonumber
\end{IEEEeqnarray}
The metric is the kinetic energy. Its matrix representation, denoted by $G$, is given by: 
\begin{equation}
G = 
\begin{pmatrix}
J_{1}+\frac{1}{4}(m_{1}+4m_{2})l_{1}^2&\frac{1}{2}m_{2}l_{1}l_{2}cos(\theta_{1}-\theta_{2})\\
\frac{1}{2}m_{2}l_{1}l_{2}cos(\theta_{1}-\theta_{2})&J_{2}+\frac{1}{4}m_{2}l_{2}^2\nonumber \\
\end{pmatrix}
\end{equation}
To write the generalized equation of geodesics, we also need the Christoffel's symbols that are given in [19] :

 $$\Gamma_{11}^{1} = \frac{m_{2}^2l_{1}^2l_{2}^2sin(2(\theta_{1}-\theta_{2}))}{8det(G)}$$

 $$\Gamma_{22}^{1} = \frac{m_{2}l_{1}l_{2}(4J_{2}+m_{2}l_{2}^2)sin(\theta_{1}-\theta_{2})}{8det(G)}$$

 $$\Gamma_{11}^{2} = - \frac{m_{2}l_{1}l_{2}(4J_{1}+(m_{1}+4m_{2})l_{1}^2)sin(\theta_{1}-\theta_{2})}{8det(G)}$$

 $$\Gamma_{22}^{2} = - \frac{m_{2}^2l_{1}^2l_{2}^2sin(2(\theta_{1}-\theta_{2}))}{8det(G)}$$
The equations of motion in absence of potential energy are :
 \begin{equation}
\ddot{\theta_{1}}+\Gamma_{11}^{1}\dot \theta_{1}^2 + \Gamma_{22}^{1}\dot \theta_{2}^2 = 0\label{eq}
\end{equation}
\begin{equation}
\ddot{\theta_{2}}+\Gamma_{11}^{2}\dot \theta_{1}^2 + \Gamma_{22}^{2}\dot \theta_{2}^2 = 0\label{eq} 
\end{equation}
\subsubsection{the tool function and it's singularities}
the tool function $x:T^{2}\rightarrow{\mathbb{R}^{2}}$ is given in local coordinates by $$x(\theta_{1},\theta_{2})=\begin{bmatrix}l_{1}cos(\theta_{1})+l_{2}cos(\theta_{2})\\l_{1}sin(\theta_{1})+l_{2}sin(\theta_{2})\end{bmatrix}$$ clearly the workspace is the ring $W_{s}=\left\{x\in\mathbb{R}^{2} / |l_{1}-l_{2}|\leq||x||_{2}\leq l_{1}+l_{2}\right\}$.\\ we compute $Dx(\theta_{1},\theta_{2})=\begin{bmatrix}-l_{1}sin(\theta_{1})&-l_{2}sin(\theta_{2})\\l_{1}cos(\theta_{1})&l_{2}cos(\theta_{2})
\end{bmatrix}$
and this gives $|det(Dx)|=l_{1}l_{2}|sin(\theta_{2}-\theta_{1})|$, the singular points are the images by the covering map $\pi:\mathbb{R}^{2}\rightarrow{T^{2}}\subset\mathbb{C}^{2}$\\ $\pi(x,y)=\begin{bmatrix}e^{ix}\\e^{iy}\end{bmatrix}$ of the set $\left\{\ x=y\right\}\cup\left\{x=y+\pi\right\}$, the singular values are the image by the tool function of the singular points, the singular values are the two circles of rays $l_{1}+l_{2}$ and $|l_{2}-l_{1}|$, for our case $l_{1}=l_{2}$ and so the singular values are the circle with ray $2l_{1}$, and the origin. 
\subsection{Application to the two link manipulators}
\subsubsection{Regulation of the tool}
We now apply this method on the two-link manipulator. Our control law is the combination of a friction force $-k.v$ and the gradient with respect to the metric of the Lasalle's potential $V(\theta_{1},\theta_{2})$.\\ \\ The control law $u$ isthusgiven by :
$$u((\theta _{1}, \theta_{2}), (\dot \theta_{1}, \dot \theta_{2})) = - grad_{g}(V(\theta_{1}, \theta_{2}))-k.v$$
Where : 
$$ V(\theta_{1},\theta_{2}) = \frac{k_{1}}{2}((l_{1}C_{1}+l_{2}C_{2}-x_{d1})^2+(l_{1}S_{1}+l_{2}S_{2}-x_{d2})^2) $$
Where $x_{d1}$ and $x_{d2}$ are respectively the first and the second components of the reference $x_{d}$.\\
\textbf{Simulation Examples :}\\
\textbf{Simulation 1} :  \\ \\  - \textit{Reference position} :  $x_{d}=(0,0.6)$   \\ - \textit{Initial conditions} : $(\theta_{1}, \theta_{2},\dot \theta_{1}, \dot \theta_{2})=(0,0,0,0)$  \\ - \textit{Gains values} : $k_{1} = 200, k =  30$  \\ \\ The results are shown in the following figures : 
\begin{figure}[H]
\centerline{\includegraphics[width=90mm]{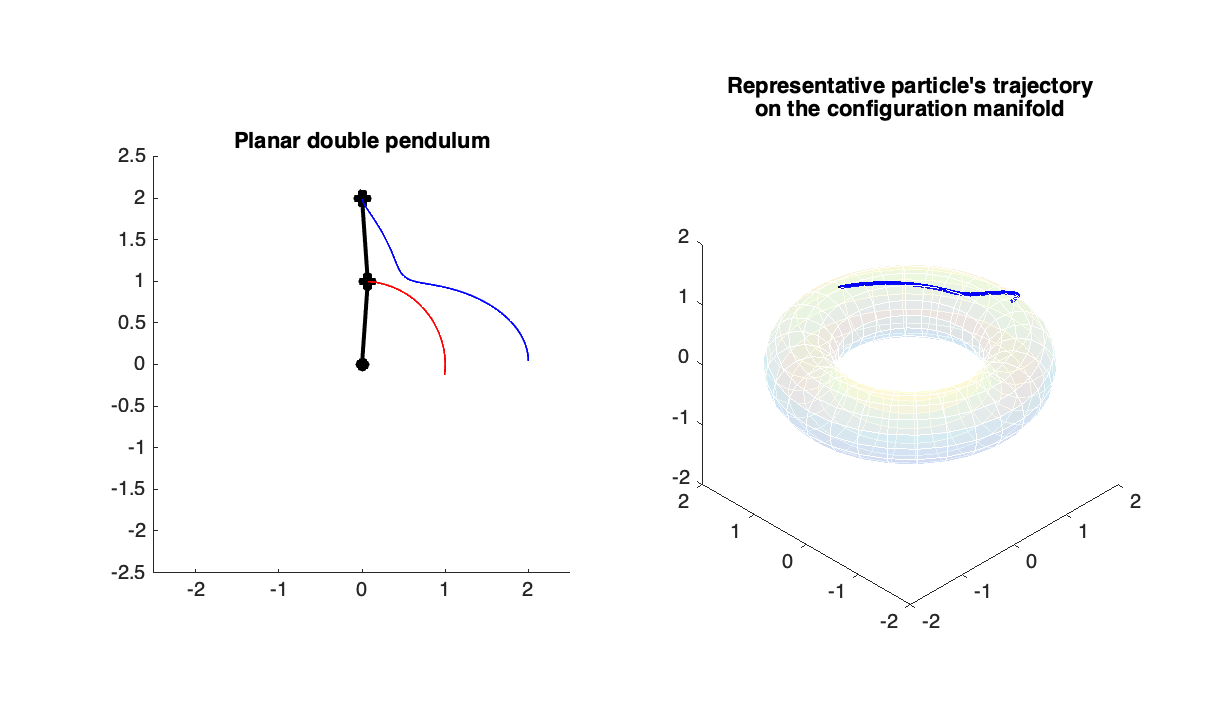}}
\caption{The two-link manipulator and the torus.}
\label{fig}
\end{figure}    

\begin{figure}[H]
\centerline{\includegraphics[width=85mm]{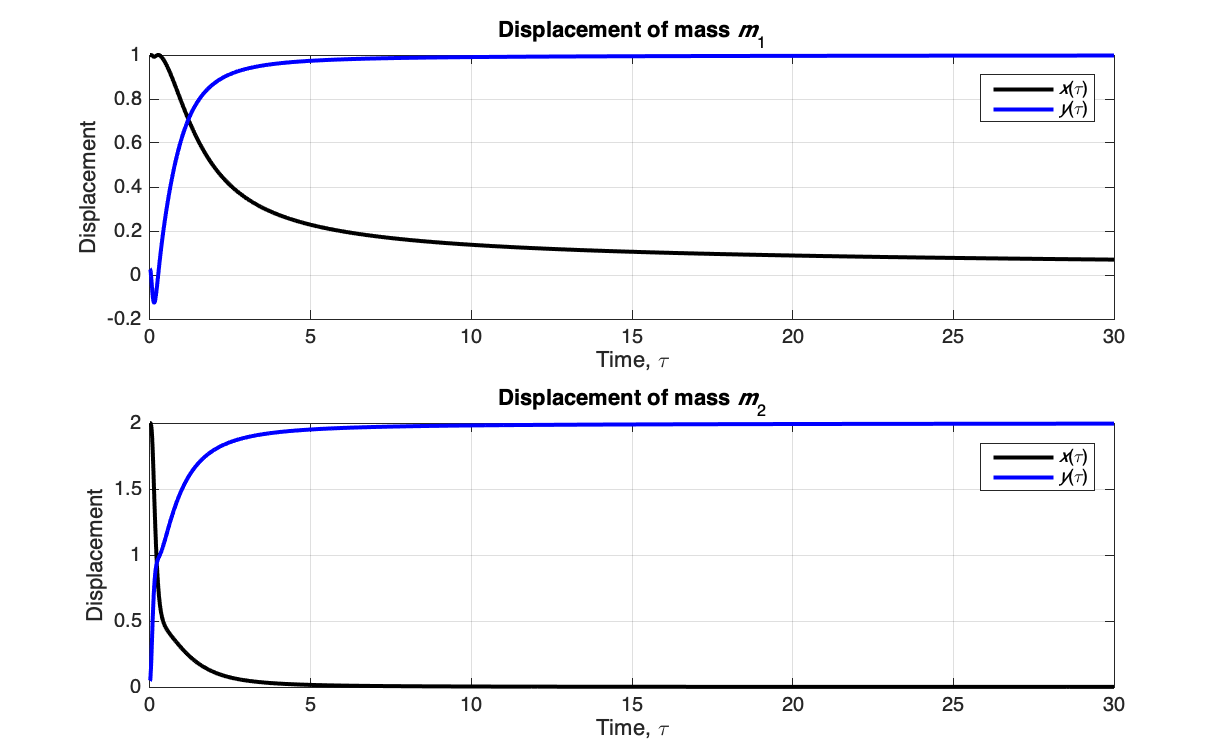}}
\caption{Trajectory of the arms.}
\label{fig}
\end{figure} 
\textbf{Simulation 2} :  \\ \\  - \textit{Reference position} :  $x_{d}=(-0.6,0)$   \\ - \textit{Initial conditions} : $(\theta_{1}, \theta_{2},\dot \theta_{1}, \dot \theta_{2})=(0,0,0,0)$  \\ - \textit{Gains values} : $k_{1} = 200, k =  30$  \\ \\ The results are shown in the following figures : 
\begin{figure}[H]
\centerline{\includegraphics[width=85mm]{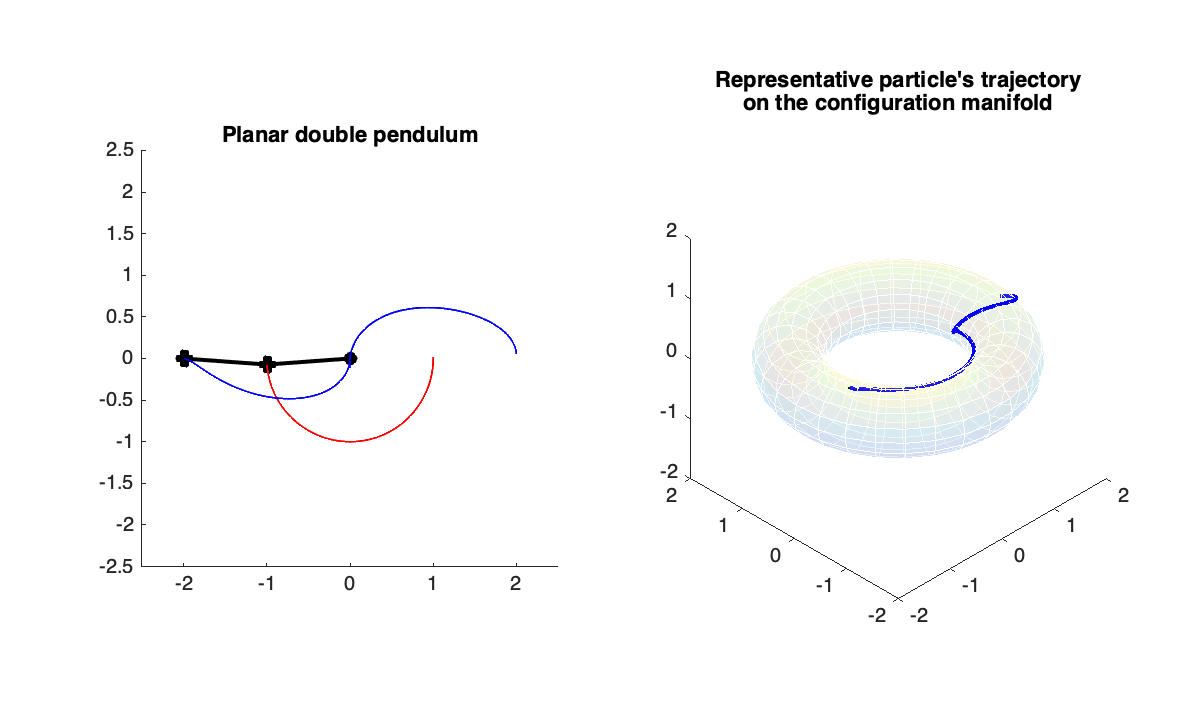}}
\caption{The two-link manipulator and the torus.}
\label{fig}
\end{figure}    
\begin{figure}[H]
\centerline{\includegraphics[width=85mm]{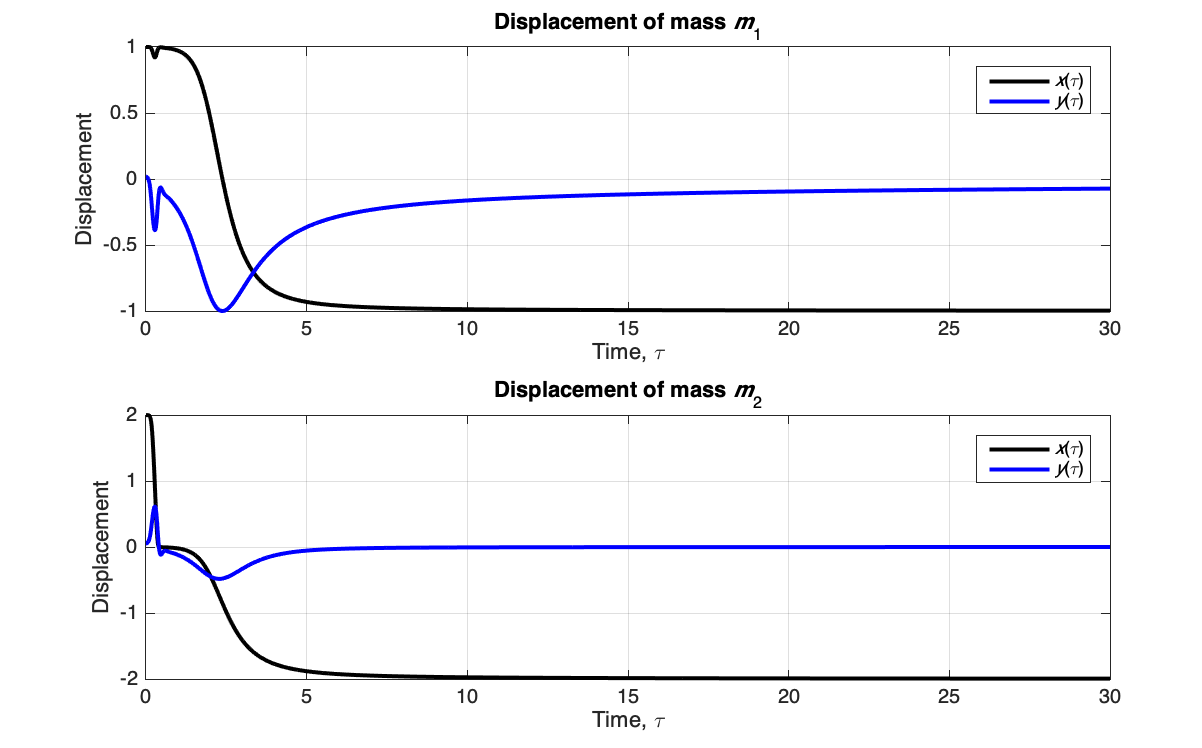}}
\caption{Trajectory of the arms.}
\label{fig}
\end{figure} 

\subsection{Constraint regulation of the tool}
We apply the developped regulator of the tool that keeps it constrained. First of all, we need to compute $\lambda$:
 $$ \lambda=\frac{-g_{\gamma}(\frac{Dgrad_{g}(\Psi)(\gamma)}{Dt},\gamma')}{|grad_{g}(\Psi)(\gamma)|_{g}^{2}} $$
 We want the surface to be a ellipse of radius $0.3$ and $0.6$ whose center is the point $(0,0)$. So $\Psi $ is given by the formula: $$\Psi (\theta_{1},\theta_{2}) = (\frac{l_{1}C_{1}+l_{2}C_{2}}{0.3})^2 + (\frac{l_{1}S_{1}+l_{2}S_{2}}{0.6})^2 - 1$$ We now calculate the gradient with respect to the metric of $\Psi $, and then we determine the covariant derivative of this gradient. The $k^{th}$ component of the covariant derivative of $grad_{g}(\Psi)$ is given by [9]: 
 \begin{IEEEeqnarray}{rCl}
\lbrace \frac{Dgrad_{g}(\Psi)(\gamma)}{Dt}\rbrace_{k}  & = & \sum_{i=1}^{2}\lbrace \frac{\partial \lbrace grad_{g}(\Psi )\rbrace_{k}}{\partial x_{i}} v_{i} \nonumber \\  
&&+ \sum_{j=1}^{2} v_{i}\Gamma _{ij}^{k} \lbrace grad_{g}(\Psi )\rbrace_{j}\rbrace \nonumber
\end{IEEEeqnarray} 

We also need to compute the norm of $grad_{g}(\Psi )$ with respect to the metric. Next, we need to compute the orthogonal projection of $grad_{g}(U_{l})$ and $v$, where $V$ is the fictive potential $$grad_{g}(V)_{//} = grad_{g}(V) - g_{q}(grad_{g}(V),n)n$$ and $$v_{//} = v - g_{q}(v,n)n$$ Where $$n=\frac{grad_{g}(\Psi )}{|grad_{g}(\Psi)|_{g}}$$ The final formula for the control law is : 
 $$u(\theta_{1},\theta_{2},\dot \theta_{1}, \dot \theta_{2}) = \lambda.grad_{g}(\Psi )-grad_{g}(V)_{//}-kv_{//}$$ 
 We can see that the norm of $grad_{g}(\Psi )$ can be zero with the presence of singularities in the torus, and we need to avoid a zeroin the denominators of all the simulations. To achieve this goal, we use some perturbation theory adding small terms $\epsilon_{1}$ and $\epsilon_{2}$ in the denominators as follow : 
 $$\lambda=\frac{-g_{\gamma}(\frac{Dgrad_{g}(\Psi)(\gamma)}{Dt},\gamma')}{|grad_{g}(\Psi)(\gamma)|_{g}^{2}+\epsilon_{1}}$$
 $$n=\frac{grad_{g}(\Psi )}{|grad_{g}(\Psi )|_{g}+\epsilon_{2}}$$
 The more $\epsilon_{1}$and $\epsilon_{2}$ are small, the more the circle is perfect.\\
 \textbf{Simulations Examples :}\\
 In the simulations we take $\epsilon_{1}=\epsilon_{2}=10^{-28}$ \\ \\ - \textit{Reference position} :  $x_{d} = (0,0.3)$   \\ - \textit{Initial conditions} : $(\theta_{1}, \theta_{2},\dot \theta_{1}, \dot \theta_{2})=(0,0,0,0)$  \\ - \textit{Gains values} : $k_{1} = 40, k =  30$ \\ \\ The results are shown in the following figures : 
\begin{figure}[H]
\centerline{\includegraphics[width=85mm]{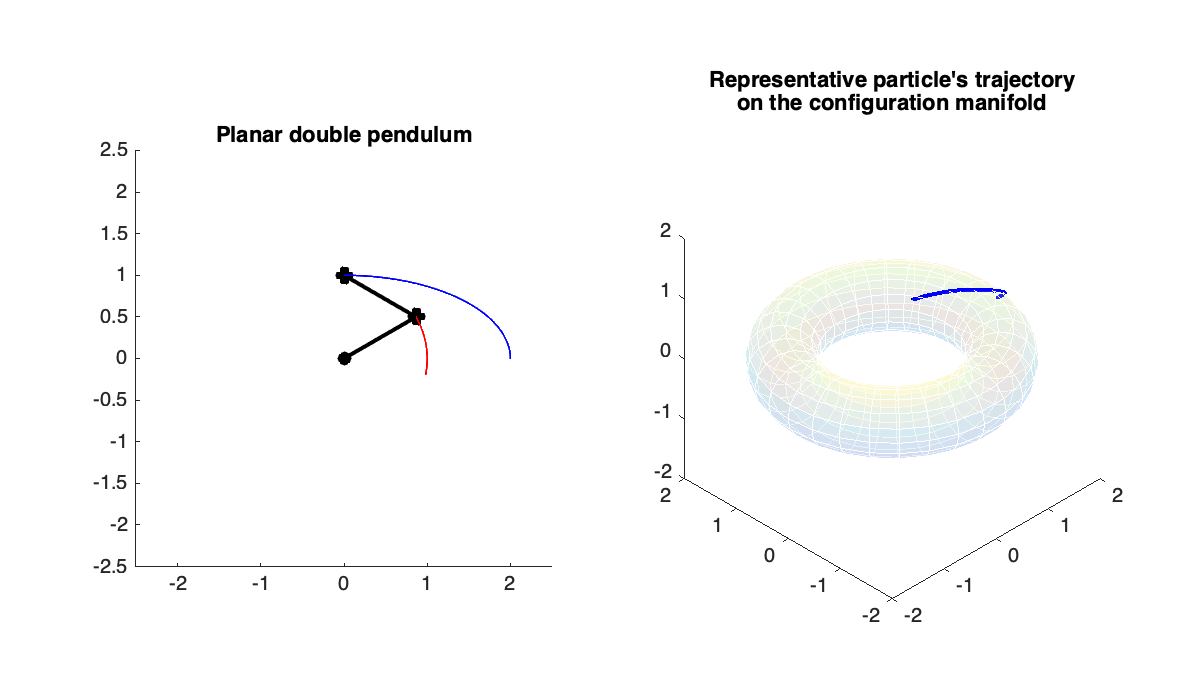}}
\caption{The two-link manipulator and the torus.}
\label{fig}
\end{figure}    

\begin{figure}[H]
\centerline{\includegraphics[width=85mm]{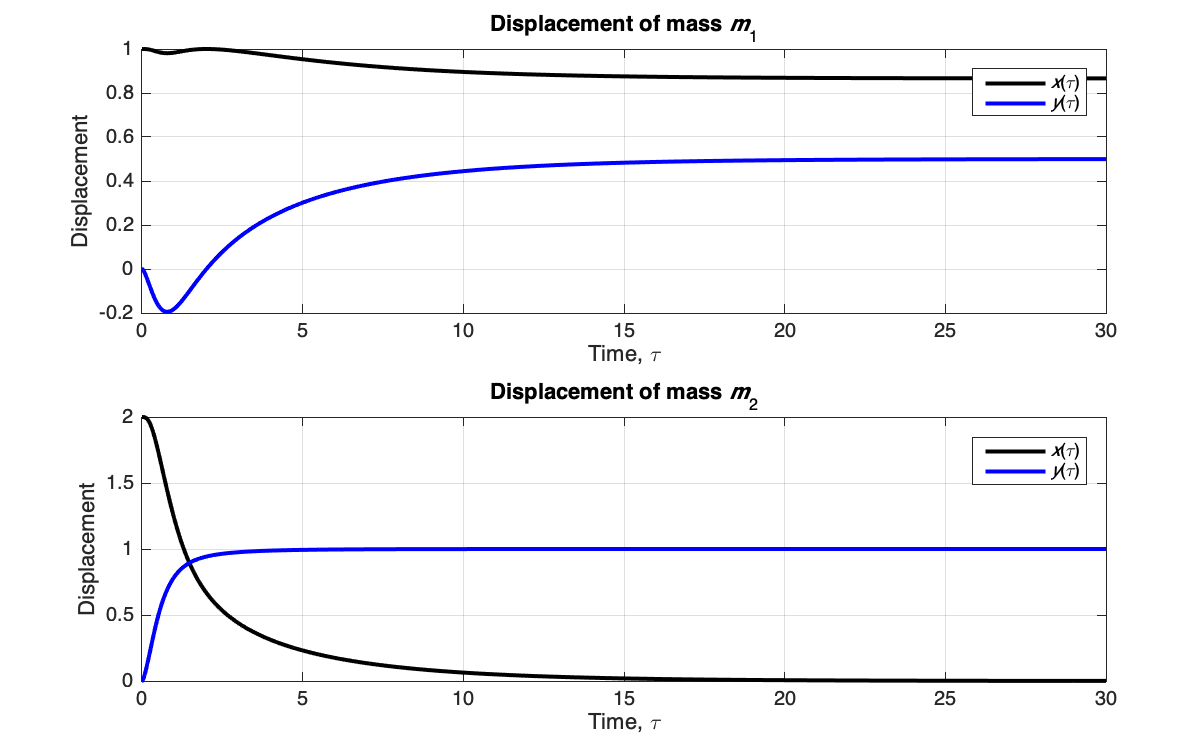}}
\caption{Trajectory of the arms.}
\label{fig}
\end{figure} 

\section{Conclusion}
Using the geometric formulation of robotic systems, instead of writing enormous equations as in the Euclidean case, we simply write basic equations that encode all the information we need, provide us with a deeper understanding of the dynamics, allow us to have a large set of control methods, and give us an exact model for our robot contrary to the Euclidean formulation. \\ \\
We give a rigorous and intrinsic formulation and proof of the tool's regulator (theorem 3.2), and for the feedback control law ensuring the geometrical constraints (theorem 3.3, theorem 3.4).\\ \\
a suite of this work can try to extend these results into robot with non-holonomic constraints as in [5] [19], or for partially actuated systems [19], and try to find some simple conditions ensuring existence and uniqueness of time and criterion optimal control, and try to give simple necessary conditions as maximum principle which allow us to compute the optimal control in practice.

\end{document}